\documentclass[reqno,11pt]{amsart}

\usepackage[a4paper,margin=1in]{geometry}

\usepackage{appendix}
\usepackage{bbm}
\usepackage{booktabs}
\usepackage{xparse}
\usepackage{sidecap}
\usepackage{float}
\usepackage{enumerate,mdwlist}
\usepackage{epsfig}
\usepackage{latexsym}
\usepackage{mleftright}
\usepackage{upgreek}

\usepackage[dvipsnames]{color}

\usepackage[bbgreekl]{mathbbol}
\usepackage{bm}

\usepackage{mathtools}
\usepackage{amssymb}
\usepackage{amsfonts}
\usepackage{amsthm}
\usepackage{aliascnt}
\usepackage{mathrsfs}
\usepackage{enumerate}
\usepackage{amscd}
\usepackage{array}
\usepackage{empheq}
\usepackage{physics}

\usepackage{subcaption}
\usepackage{tikz,tikz-cd}
\usepackage{tkz-euclide}
\usetikzlibrary{matrix}
\usepackage{smartdiagram}
\usepackage{graphicx}
\usepackage{longtable}
\usepackage{multirow}

\usepackage{etoolbox}
\usepackage[export]{adjustbox}

\numberwithin{equation}{section}

\usepackage{hyperref} 
\hypersetup{
	colorlinks = true,
	linkcolor = blue,
	anchorcolor = blue,
	citecolor = blue,
	filecolor = blue,
	urlcolor = blue
}
\usepackage[utf8]{inputenc}
\usepackage{url}
\usepackage[backend=biber,style=alphabetic,maxnames=99,minnames=99]{biblatex}
\usepackage[capitalize,noabbrev,nameinlink]{cleveref}

\renewbibmacro{in:}{%
	\ifentrytype{article}
	{}
	{\printtext{\bibstring{in}\intitlepunct}}}

\makeatletter
\patchcmd\maketitle
{\uppercasenonmath\shorttitle}
{}
{}{}
\patchcmd\maketitle
{\@nx\MakeUppercase{\the\toks@}}
{\the\toks@}
{}{}
\patchcmd\@settitle
{\uppercasenonmath\@title}
{}
{}{}
\patchcmd\@setauthors
{\MakeUppercase{\authors}}
{\authors}
{}{}
\makeatother

\crefname{theorem}{Theorem}{Theorems}
\crefname{definition}{Definition}{Definitions}
\crefname{fact}{Fact}{Facts}
\crefname{proposition}{Proposition}{Propositions}
\crefname{lemma}{Lemma}{Lemmas}
\crefname{corolary}{Corolary}{Corolaries}
\crefname{assumption}{Assumption}{Assumptions}
\crefname{claim}{Claim}{Claims}
\crefname{result}{Result}{Results}
\crefname{problem}{Problem}{Problems}
\crefname{remark}{Remark}{Remarks}
\crefname{example}{Example}{Examples}
\crefname{corollary}{Corollary}{Corollaries}
\crefname{subsubsection}{Subsubsection}{Subsubsections}
\crefname{subsection}{Subsection}{Subsections}
\crefname{section}{Section}{Sections}

\newcommand{\bbR}{\mathbb{R}}
\newcommand{\bbE}{\mathbb{E}}

\newcommand{\bbN}{\mathbb{N}}
\newcommand{\bbT}{\mathbb{T}}
\newcommand{\bbZ}{\mathbb{Z}}

\newcommand{\calP}{\mathcal{P}}

\newcommand{\calG}{\mathcal{G}}

\newcommand{\calM}{\mathcal{M}}

\newcommand{\calR}{\mathcal{R}}
\newcommand{\calL}{\mathcal{L}}

\newcommand{\bfX}{\mathbf{X}}

\newcommand{\bfP}{\mathbf{P}}
\newcommand{\bfQ}{\mathbf{Q}}
\newcommand{\bfD}{\mathbf{D}}

\newcommand{\frakP}{\mathfrak{P}}

\newcommand{\scrB}{\mathscr{B}}

\newcommand{\scrF}{\mathscr{F}}
\newcommand{\sfd}{\mathsf{d}}

\newcommand{\Vol}{\operatorname{vol}}
\newcommand{\Isom}{\operatorname{Isom}}

\newcommand{\Stab}{\operatorname{Stab}}

\newcommand{\Id}{\operatorname{Id}}

\newcommand{\LCcnt}{\nabla^{\operatorname{LC}}}

\newcommand{\gfr}{g_{\operatorname{FR}}}
\newcommand{\gqout}{g_{\operatorname{quot}}}
\newcommand{\gqfr}{g_{\operatorname{inv}}}
\newcommand{\Var}{\operatorname{Var}}
\newcommand{\cac}{C_{\operatorname{AC}}}
\newcommand{\AlphaC}{\nabla^{(\alpha)}}

\theoremstyle{plain}
\newtheorem{theorem}{Theorem}[section]

\newaliascnt{proposition}{theorem}
\newtheorem{proposition}[proposition]{Proposition}
\aliascntresetthe{proposition}

\newaliascnt{lemma}{theorem}
\newtheorem{lemma}[lemma]{Lemma}
\aliascntresetthe{lemma}

\newaliascnt{corollary}{theorem}
\newtheorem{corollary}[corollary]{Corollary}
\aliascntresetthe{corollary}

\newaliascnt{claim}{theorem}

\aliascntresetthe{claim}

\theoremstyle{definition}
\newaliascnt{definition}{theorem}
\newtheorem{definition}[definition]{Definition}
\aliascntresetthe{definition}

\newaliascnt{problem}{theorem}

\aliascntresetthe{problem}

\theoremstyle{remark}

\newaliascnt{remark}{theorem}
\newtheorem{remark}[remark]{Remark}
\aliascntresetthe{remark}

\newaliascnt{example}{theorem}
\newtheorem{example}[example]{Example}
\aliascntresetthe{example}

\begin{document}


	\title[Group Invariance and Isometry Subgroup]{Group Invariance of Statistical Models and Associated Statistical Manifolds}
	
	\author[Z. Cheng]{Zhuohao Cheng}
	
	\address{Graduate School of Mathematics \endgraf
		Nagoya University \endgraf
		464-0812, Chikusa, Nagoya, Aichi, Japan}
	
	\email{cheng.zhuohao.y3@s.mail.nagoya-u.ac.jp}

	\begin{abstract}
		We study the group invariance of a parameterized statistical model that is closed under push-forward sample transformations given by a group action. 
		We prove that the group action on the sample space induces an isometric action on the associated statistical manifold. Furthermore, it preserves the Fisher-Rao metric as well as the Amari-Chentsov tensor. 
		Moreover, we establish that this isometric action is compatible with equivariant and sufficient statistics in the spirit of the Hall-Wijsman-Ghosh generalized Stein's theorem. 
		This compatibility leads to a novel reduction of the original group-invariant model, where its Fisher information loss is explicitly quantified via conditional expectation. 
		As an application, we apply this framework to a generalized cardioid distribution in circular data analysis, illustrating how the reduced model facilitates estimation and hypothesis testing for relative parameters.
	\end{abstract}
	
	
	\maketitle

	\section{Introduction}\label{section_introduction}
	
	Information geometry has been an irreplaceable tool in probability theory, statistics, and beyond \cite{Amari2016,Ayetal2017,Chewi2025}.
	It provides geometric means to study statistical theories; see, e.g., \cite{Ay2015,Le2006,Pollard2013,Yamaguchi2024} for more details. 
	Recent studies have explored various extensions of statistical manifolds and the Fisher-Rao metric; see, e.g., \cite{Kainth2025,Ay2025} and references therein. 
	In applications, many parameterized statistical models naturally have some symmetries. That is, the model is invariant under the parameter transformations induced by the group action on the sample space. 
	Recently, such group-invariant statistical models, manifolds, and related subjects have been studied from various perspectives \cite{Neilson2026,Inoguchi2025,Coatanhay2026,Kobayashi2026};
	Inoguchi and Ohno \cite{Inoguchi2025} give a comprehensive study on \textit{homogeneous statistical manifolds}. They develop a general theory of homogeneous manifolds from the viewpoint of Fisher-Rao geometry.
	Kobayashi et al. \cite{Kobayashi2026} build the theory of moduli spaces of left-invariant statistical structures on Lie groups that is the central concept in  \cite{Furuhata2021}. These studies provide theoretical frameworks for analyzing symmetric properties of statistical manifolds. 
	Particularly, Neilson and Okamura \cite{Neilson2026} study the statistical model that is equivariant under group action.
	They assume that a given group can act on the sample space and also on the parameter space, and the density functions transform with appropriate multipliers.
	Then, they study the invariance of $f$-divergences that generate the associated statistical manifold, and find that such invariance of the model yields a reduction to the divergence manifold.
	
	Inspired by these studies, in this paper, we investigate the geometric structures of the statistical manifold associated with a group-invariant statistical model. 
	Specifically, we link the group invariance of the model to an isometry subgroup of the associated statistical manifold. 
	We find that this isometry subgroup always forms a Lie group acting on the manifold while preserving its statistical structures, and remains invariant under equivariant and sufficient statistics. 
	Moreover, several classical results on invariant and sufficient inference, see, e.g., \cite{Hall1965,Chiu2025}, can be revisited from this geometric perspective. 
	We further show that this group invariance induces a reduction of the model and quantify the resulting information loss. 
	The detailed contributions and organization of this paper are summarized as follows.
	In \cref{section_preliminaries}, we review some necessary notions about parameterized statistical models and \textit{Hall-Wijsman-Ghosh's generalized Stein's theorem on invariance and sufficiency}; see \cite{Hall1965}.
	In \cref{subsection_geometric_interpretation_of_equivariant_sufficiency}, we first consider a geometric interpretation of the group invariance of the statistical model from the perspective of the isometry subgroup. 
	We adopt some techniques in general Banach manifolds to proceed with our arguments in terms of conditional expectation, in which the theory of immersion and submersion between Banach manifolds helps us establish the regularity of induced parameter transformation.
	This is helpful when we interpret the results of the present paper from the viewpoint of statistics.
	\cref{proposition_isometry_group_for_FR_metric}-\cref{theorem_geometric_analogy} show that there exists a group homomorphism between the transformation group acting on the sample space and the isometry group on the statistical manifold. Furthermore, the induced isometry subgroup can be preserved under equivariant and sufficient statistics.
	\cref{subsection_quotient_Wasserstein_versus_Wasserstein_of_quotient} extends the invariant property of the Fisher-Rao metric to the Amari-Chentsov tensor and hence $\alpha$-connections. 
	This extension shows that whenever the induced isometry subgroup acts transitively, the statistical manifold automatically becomes homogeneous in the sense of \cite{Inoguchi2025}.
	
	In particular, we find that the group invariance of a model yields a reduction on the original model in a distinctive way from the reduction given in \cite{Neilson2026}, in which the reduction yields a quotient space where the divergences and Fisher-Rao distance can naturally descend without contraction.
	The reduction studied in the present paper yields a new parameterized statistical model, which only depends on the invariant observation.
	This motivates us to study the information loss caused by this reduction and the relation between the reduction and (equivariant and) sufficient statistics. 
	In \cref{section_invariant_statistical_structures_and_homogeneity}, we deduce the \textit{reduced statistical model} by applying the \textit{restricted model on invariant observation}. 
	By applying some similar techniques as in the previous section, we study the reduced model from the perspective of conditional expectation.
	Then, we consider the information loss of the reduction.
	\cref{theorem_comparison_of_quotient_and_induced_FR}-\cref{theorem_estimation_of_information_loss_caused_by_restriction} quantify the horizontal and total information loss and link them to the conditional variance of the invariant observation.
	
	In the final section, together with a concrete example, we will see an application to circular data analysis. We use a generalized cardioid distribution to construct a non-trivial reduced model and calculate the closed form of the information loss in this case.

	\paragraph{\textbf{Notations and Conventions}} Throughout the paper, we will use $(X,\sfd)$ to represent a metric space equipped with the metric topology induced by the distance function $\sfd$. Let $\mathbb{N}$ be the set of natural number including $0$, and $\mathbb{N}_+=\mathbb{N}\setminus\{0\}$. We will use standard notations such as $\mathbb{R}^n$ for $n$-dimensional Euclidean space, and $(M,g)$ for a Riemannian manifold with the metric $g$ on the tangent space $T_pM$ for each $p\in M$. We use $(X,\mathscr{B})$ to present a measurable space. We use $\calP(X)$ to represent the set of probability measures on $X$. In the present paper, when we use the terminology ``\textit{sufficiency}'', it implicitly implies the classical definition, cf. \cref{definition_sufficient_sigma_algebra}.

	\section{Preliminaries on Statistical Models and Conditional Expectation}\label{section_preliminaries}

	We first recall some background knowledge about parameterized statistical models. Then we adopt some necessary notions about conditional expectation.
	Throughout the paper, we assume that the Riemannian manifold $(M,g)$ is complete, smooth, $n$-dimensional, connected, and without boundary or corners.
	For any measurable space $(X,\scrB)$, we always assume $\scrB\subset\scrB_X$, i.e., the subalgebra of Borel $\sigma$-algebra of $X$ with respect to the metric topology.

	\subsection{Parameterized statistical model}\label{subsection_parameterized_statistical_model}
	
	We adopt some standard notions in information geometry from \cite{Amari2000, Ayetal2017}.
	
	\begin{definition}\label{definition_parameterized_statistical_model}
		We say that $(M,X,\mathbf{P})$ is a \textit{parameterized statistical model} if $M$ is a general manifold, $X$ is the underlying space of a measurable space $(X,\mathscr{B})$, and $\mathbf{P}:M\to\mathcal{P}(X)$ is a $C^1$-mapping in the sense of \cite[Definition C.3]{Ayetal2017}.
	\end{definition}
	We briefly note that the theory of general differentiable maps between Banach spaces or manifolds can be found in \cite{Abraham1988,Schmeding2023} as well. The next notion that we need is the \textit{statistical manifold}.  
	\begin{definition}\label{definition_statistical_manifold}
		A general statistical manifold is the triple $(M,g,C)$ consisting of a manifold $M$, a Riemannian metric $g$ and a $3$-tensor $C$ that is symmetric in all arguments.
	\end{definition}
	In the present paper, when we use the terminology ``statistical manifold'', we mainly consider a Riemannian manifold equipped with the Fisher-Rao metric induced by the associated parameterized statistical model, cf. \cref{definition_Fisher_Rao_metric} below.
	Therefore, we will simply use the notation $(M,\gfr)$ to represent a statistical manifold equipped with this metric.
	Occasionally, we will use the notation $(M,\gfr,\cac)$ where $\cac$ is the \textit{Amari-Chentsov cubic tensor} rather than a general totally symmetric $(0,3)$-tensor.
	
	To give the definition of sufficient statistics and the Fisher-Rao metric in the context of a parameterized model, we need some additional notions. 
	We define the score function of $\bfP$ at $\theta\in M$ for any $V\in T_{\theta}M$ by  $s_{\theta,V}=\frac{d(D_{\theta}\mathbf{P}(V))}{d\mathbf{P}(\theta)}=:\partial_V\log\mathbf{P}(\theta)\in L^1(X,\mathbf{P}(\theta))$ where $D_{\theta}\mathbf{P}(V)$ is the image of a tangent vector $V\in T_{\theta}M$ via the differential $D\mathbf{P}$ of $\mathbf{P}$ at $\theta$, which is a finite signed measure with integral $0$ on $X$. See, e.g., \cite[Chapter 3.2, Chapter 5.1]{Ayetal2017} for details. 
	
	Throughout the paper, we always assume the model $(M, X,\mathbf{P})$ is $k$-integrable in the sense of \cite[Definition 3.7]{Ayetal2017} for some  $k\geq3$. In particular, $s_{\theta,V}\in L^k(X,\bfP_{\theta})$ for any $\theta\in M$, $V\in T_{\theta}M$. 
	We note that the assumption of $k$-integrability does not need to imply the smoothness of $\gfr$ (and $\cac$); see \cite{Yamaguchi2024,Inoguchi2025}. Therefore, unless we explicitly mention otherwise, we assume that the induced Fisher-Rao tensor is smooth and positive-definite, i.e., it becomes a smooth metric on $M$.
	
	\begin{definition}\label{definition_Fisher_Rao_metric}
		Let $(M,X,\mathbf{P})$ be a parameterized statistical model. Then the Fisher-Rao metric and Amari-Chentsov tensor defined on the tangent space $T_{\theta}M$ of $M$ at $\theta\in M$ are given by 
		$$
		\gfr(U,V)_{\theta}:=\int_{X}s_{\theta,U}s_{\theta,V}d\mathbf{P}_{\theta}(x),\ \cac(U,V,W)_{\theta}:=\int_{X}s_{\theta,U}s_{\theta,V}s_{\theta,W}d\mathbf{P}_{\theta}(x)
		$$
		for any $U,V,W\in T_{\theta}M$ where $\mathbf{P}_{\theta}:=\mathbf{P}(\theta)\in\mathcal{P}(X)$.
	\end{definition}
	
	\begin{definition}\label{definition_geometric_sufficiency}
		Let $(M,X,\mathbf{P})$ be a parameterized statistical model. Then, a statistic $T: X\to X$ is called a (geometric) sufficient statistic if for any tangent vector $V\in T_{\theta}M$ there is
		$$
		||\partial_V\log T_{\#}\mathbf{P}(\theta)||_{L^k(X,T_{\#}\bfP_{\theta})}=||\partial_V\log\mathbf{P}(\theta)||_{L^k(X,\bfP_{\theta})}
		$$
		for any $\theta\in M$.
	\end{definition}
	
	\subsection{Preliminaries on transformations and conditional expectation}\label{subsection_preliminaries_on_transformations_and_conditional_independence}
	Let us recall some preliminaries adopted from \cite{Hall1965}\footnote{There is a correction on the paper provided by the authors; see \cite{Hall1995}. The correction concerns the assumption of Lemma 3.3 in the original paper, while we will not use it and related results.}. For a rigorous theory in statistics, see, e.g., \cite{Lehmann2005,Dryden2016}. For the detailed theory on conditional expectation and measure theory used below, see, e.g., \cite[Chapter 10, Vol.II]{Bogachev2007}.
	
	Let $T: X\to X$ be a Borel map (i.e., a statistic), $G$ a group of sample transformations on $X$ (for tractability). We say that $T$ is an \textit{equivariant statistic} with respect to $G$ if for every $h\in G, x\in X$, there is $T(hx)=hT(x)$, i.e., the group action and statistic commute. 
	
	Denote a measurable space and a family of probability measures on $X$ by $(X,\mathscr{B})$ and $\frakP\subset\calP(X)$, respectively. We say that $\mathscr{B}_0$ is a subalgebra of $\mathscr{B}$ if $\mathscr{B}_0\subset\mathscr{B}$ and $\mathscr{B}_0$ is a $\sigma$-algebra. 
	For two $\mathscr{B}$-measurable functions $f_1,f_2:X\to\mathbb{R}$, we say that $f_1$ is \textit{almost equivalent} to $f_2$, denoted by $f_1\sim f_2$, with respect to some $\mu\in\frakP$ if $f_1=f_2$ $\mu$-a.e. on $X$. Similarly, we say that $f_1\sim f_2$ with respect to $\frakP$ if $f_1=f_2$ $\mu$-a.e. for all $\mu\in\frakP$. 
	
	\begin{definition}\label{definition_conditional_expectation}
		Let $\mu$ be a probability measure on $(X,\mathscr B)$, let $\mathscr B_0\subseteq\mathscr B$ be a subalgebra, and let $f\in L^1(X,\mathscr B,\mu)$. The \textit{conditional expectation} of $f$ given $\mathscr B_0$ is a $\mathscr B_0$-measurable integrable function $f_0$
		such that
		$$
		\int_A f_0\,d\mu=\int_A f\,d\mu,\ \ A\in\mathscr B_0.
		$$
		It is unique up to $\mu$-a.e. equality and is denoted by $\mathbb E_\mu(f\mid\mathscr B_0)$.
	\end{definition}

	Denote the set of bijective maps on $X$ by $\operatorname{Inv}(X)$. We therefore assume that $G\leq\operatorname{Inv}(X)$ is a group under map composition.
	The point transformation $h$ induces naturally an action on $\mathcal{P}(X)$ (and hence on $\frakP\subset\mathcal{P}(X)$) by push-forward of measures. 
	For the measurable space $(X,\mathscr{B})$, define the $\sigma$-algebra $h\mathscr{B}$ by $B\in h\mathscr{B}\iff h^{-1}(B)\in\mathscr{B}$. 
	
	\begin{definition}\label{definition_sufficient_sigma_algebra}
		For any two subalgebras $\mathscr{B}_1$ and $\mathscr{B}_2$ such that $\mathscr{B}_2\subset\mathscr{B}_1$, we say that $\mathscr{B}_2$ is sufficient for $\mathscr{B}_1$ if for every bounded $\mathscr{B}_1$-measurable function $f_1$, there is a bounded $\mathscr{B}_2$-measurable $f_2$ such that $\mathbb{E}_{\mu}(f_1|\mathscr{B}_2)\sim f_2$ $\mu$-a.e. for any $\mu\in\frakP$. In particular, if a Borel subalgebra $\mathscr{B}_S\subset\mathscr{B}$ is sufficient, and $T:X\to X$ is the map such that $\scrB_S=T^{-1}(\scrB)$, then $T$ is (classical) sufficient.
	\end{definition}
	
	\begin{definition}\label{definition_invariant_sigma_algebra}
		We say that a subset $A\in\mathscr{B}$ of $X$ is invariant under $G$ if its indicator function is invariant, i.e., $hA=A$. Furthermore, we say that the $\sigma$-algebra $\mathscr{B}_I$, generated by the collection of invariant sets, is the invariant $\sigma$-algebra of $X$ under $G$.
	\end{definition}
	
	To introduce the main result in \cite{Hall1965}, let us first recall some necessary notions in that paper. For a measurable space $(X,\scrB)$ and a group action $G\leq\operatorname{Inv}(X)$, we will always assume that $h\in G$ is Borel measurable throughout the paper.
	
	\begin{definition}\label{definition_group_invariant_model}
		Let $\frakP\subset\calP(X)$ be a family of probability measures on the measurable space $(X,\scrB)$, $G\leq\operatorname{Inv}(X)$. We say that $\frakP$ is $G$-invariant if $h_{\#}\frakP=\frakP$ for any $h\in G$. Here, $h_{\#}\frakP:=\{ h_{\#}\mu\ |\ \mu\in\frakP \}$ and $h_{\#}\mu$ is the push-forward of $\mu$ via $h$, see, e.g., \cite[Chapter 1-2]{Figalli2023}.
	\end{definition}
	
	We say that a measurable function $f$ is \textit{almost invariant} if, for every $h\in G$ and $\mu\in P$,
	$$
	f\circ h=f\ \ \mu\operatorname{-a.e.}
	$$
	Throughout the section, we will define the subalgebra $\mathscr{B}_{SI}=\mathscr{B}_S\cap\mathscr{B}_I$ where $\mathscr{B}_I$ is the invariant subalgebra of $\mathscr{B}$ under $G$ and $\mathscr{B}_S$ is the sufficient subalgebra of $\mathscr{B}$ induced by a sufficient statistic $T:X\to X$.
	
	\begin{theorem}[{\cite[Theorem 3.1]{Hall1965}}]\label{ref.Hall_theorem_3.1}
		Let $(X,\scrB)$ be a measurable space, $G\leq\operatorname{Inv}(X)$, and $\frakP$ a $G$-invariant family. 
		If $h\scrB_S=\scrB_S$ for any $h\in G$, and for any almost invariant $\scrB_S$-measurable $f_S$ with respect to $\frakP$, there exists a $\scrB_{SI}$-measurable function $f_{SI}$ such that $f_{SI}\sim f_{S}$ with respect to $\frakP$.
		Then $\mathscr{B}_{SI}\subset\mathscr{B}_I$ is sufficient.
	\end{theorem}
	
	The preceding theorem explains the following commutative diagram; see \cite[Part I.1]{Hall1965} for details.
	\[
	\begin{tikzcd}[column sep=huge,row sep=large]
		\bfX_{\theta}
		\arrow[r]
		\arrow[d]
		&
		\mathbf{S}_{\theta}
		\arrow[d]
		\\
		\mathbf{U}_{\gamma}
		\arrow[r]
		&
		\mathbf{V}_{\gamma}.
	\end{tikzcd}
	\]
	Here, the vertical and horizontal arrows mean the maximal invariant reduction and sufficient reduction, respectively.
	$\bfX_{\theta}$ may be thought of as the random variable determined by $\theta$; $\mathbf{S}_{\theta}$ is the random variable of a sufficient statistic defined by $\bfX_{\theta}$. 
	$\mathbf{U}_{\gamma}$ may be thought of as the probability model for a maximal invariant function $u$, which is defined by $\bfX_{\theta}$ and only depends on $\gamma$; hence $\mathbf{V}_{\gamma}$ is the sufficient reduction of $\mathbf{U}_{\gamma}$.


	\section{Isometries Induced by Invariant Reduction}\label{section_invariance_sufficiency_and_isometry}
	
	In this section, we first investigate the relation between group invariance of statistical models, cf. \cref{definition_group_invariant_model}, and isometries on the associated statistical manifolds. 
	We will first show a compatibility between the induced isometry subgroup and equivariant and sufficient statistics in the spirit of Hall-Wijsman-Ghosh's theorem; then we will see that the induced closed isometry subgroup preserves the statistical structures on the manifold.
	
	A necessary notion for this section is the \textit{$G$-invariant statistical model}.
	\begin{definition}\label{definition_G_invariant_statistical_model}
		We say that a parameterized statistical model $(M,X,\bfP)$ is $G$-invariant with respect to $G\leq\operatorname{Inv}(X)$ having Borel measurable elements if it holds
		$$
		h_{\#}\bfP(M):=\{ h_{\#}\bfP(\theta)\ |\ \theta\in M \}=\bfP(M)
		$$
		for any $h\in G$. Here $h_{\#}\bfP(\theta)$ is the push-forward of $\bfP(\theta)$ via the Borel map $h$.
	\end{definition}
	
	We will first give our general assumption.
	In the remainder of the paper, for some positive integer $k\geq3$, we assume that any $k$-integral model $(M,X,\bfP)$ under consideration is $G$-invariant with respect to $G$, where $G\leq\operatorname{Inv}(X)$ is a group action on the associated sample space $(X,\scrB)$ with Borel measurable elements.
	Furthermore, we assume that the model is \textit{identifiable}, i.e., $\theta\mapsto\bfP(\theta)$ is injective.

	\subsection{Group invariance and isometry subgroup}\label{subsection_geometric_interpretation_of_equivariant_sufficiency}
	
	Let us first investigate the relation between group invariance of the statistical model and isometries on the associated statistical manifold. 
	
	We first investigate some properties of the parameter transform induced by $h\in G\leq\operatorname{Inv}(X)$. We denote the space of \textit{finite signed measures} on a measurable space $X$ by $\calM(X)$. This is a Banach space equipped with the total variation norm.
	
	\begin{lemma}\label{lemma_the_transport_of_Radon_Nikodym_derivative}
		Let \(h:X\to X\) be a measurable bijection with measurable inverse, and let \(\mu\) and \(\nu\) be \(\sigma\)-finite measures in $\calM(X)$. If \(\mu\ll\nu\), then \(h_{\#}\mu\ll h_{\#}\nu\) and
		$$
		\frac{d(h_{\#}\mu)}{d(h_{\#}\nu)}=\frac{d\mu}{d\nu}\circ h^{-1}\ \ h_{\#}\nu\text{-a.e.}
		$$
	\end{lemma}
	
	\begin{proof}
		Write \(p=d\mu/d\nu\). For every bounded Borel function \(\varphi\),
		\begin{align*}
			\int_X \varphi(x)\,d(h_{\#}\mu)(x) & =\int_X \varphi\circ h(x)\,d\mu(x)=\int_X (\varphi\circ h)(x)\cdot p(x)\,d\nu(x) \\
			& =\int_X \varphi(x)\cdot\,(p\circ h^{-1})(x)\,d(h_{\#}\nu)(x).
		\end{align*}
		The claimed identity follows from the uniqueness of the Radon-Nikodym derivative; see, e.g., \cite[Theorem 6.10]{Bercovici2016}. This shows $h_{\#}\mu\ll h_{\#}\nu$.
	\end{proof}
	
	\begin{proposition}\label{proposition_isometry_group_for_FR_metric}
		Let $(M,X,\bfP)$ be a parameterized statistical model and $G\leq\operatorname{Inv}(X)$.
		For any $h\in G$, there exists a parameter transformation $\varphi_h:M\to M$ such that $h_{\#}\bfP(\theta)=\bfP(\varphi_h(\theta))$ for any $\theta\in M$ and $\varphi_h\in\Isom(M,\gfr)$.
	\end{proposition}
	
	\begin{proof}
		We denote the score function of $\bfP_\theta$ at $\theta\in M$ with respect to $V\in T_\theta M$ by
		$$
		s_{\theta,V}=\frac{d(D_\theta\bfP(V))}{d\bfP(\theta)}=:\partial_V\log\bfP(\theta).
		$$
		By the regularity assumptions, $\bfP$ is a $C^1$ map into $\calM(X)$ equipped with the total variation norm.
		
		Since $\gfr$ is positive-definite, we have
		$$
		D_\theta\bfP(V)=0 \ \implies\ s_{\theta,V}=0\ \bfP_\theta\operatorname{-a.e.}\ \implies\ \gfr(V,V)_\theta=0 \ \implies\ V=0.
		$$
		Hence $\operatorname{Ker}D_\theta\bfP=\{0\}$. Since $M$ is finite-dimensional, $\bfP$ is a $C^1$-immersion.
		
		Next, for any $h\in G$, recall that both $h,h^{-1}:X\to X$ are Borel. We can therefore define a bounded linear operator $L_h:\calM(X)\to\calM(X)$ by
		$$
		L_h(\mu):=h_\#\mu.
		$$
		Its inverse is $L_{h^{-1}}$. Since $M$ is compact and $\calM(X)$ is Hausdorff, the continuity and injectivity of $\bfP$ imply that $\bfP:M\to\bfP(M)$ is a homeomorphism onto its image.
		Combining this with the immersion property, we conclude that $\bfP$ is a $C^1$-embedding. Thus $\bfP(M)$ is a $C^1$ embedded submanifold of $\calM(X)$, and
		$$
		\bfP^{-1}:\bfP(M)\to M
		$$
		is $C^1$. The $C^1$ regularity of $\bfP^{-1}$ follows from \cite[Theorem H, Lemma 1.13]{Glockner2015} or \cite[Chapter 1.51, 1.56]{Schmeding2023}.
		
		By the $G$-invariance of the model, $L_h(\bfP(M))=\bfP(M)$. Hence we can define
		$$
		\varphi_h:=\bfP^{-1}\circ L_h\circ\bfP:M\to M.
		$$
		This map is $C^1$, and its inverse is $\varphi_h^{-1}=\bfP^{-1}\circ L_{h^{-1}}\circ\bfP=\varphi_{h^{-1}}$, which is also $C^1$. Therefore, $\varphi_h$ is a $C^1$-diffeomorphism.
		
		Differentiating the identity $\bfP\circ\varphi_h=L_h\circ\bfP$, we obtain
		$$
		D_{\varphi_h(\theta)}\bfP\bigl(D_\theta\varphi_h(V)\bigr)=h_\#\bigl(D_\theta\bfP(V)\bigr).
		$$
		Applying \cref{lemma_the_transport_of_Radon_Nikodym_derivative} to the finite signed measure $D_\theta\bfP(V)$ gives
		$$
		s_{\varphi_h(\theta),D_\theta\varphi_h(V)}=\frac{d(h_{\#}(D_{\theta}\bfP(V)))}{d(h_{\#}\bfP_{\theta})}=s_{\theta,V}\circ h^{-1}\ \ \bfP_{\varphi_h(\theta)}\operatorname{-a.e.}
		$$
		Since $\bfP_{\varphi_h(\theta)}=h_\#\bfP_\theta$, it follows that
		\begin{align*}
			\gfr\bigl(D_\theta\varphi_h(V),D_\theta\varphi_h(W)\bigr)_{\varphi_h(\theta)} & =\int_X(s_{\theta,V}s_{\theta,W})\circ h^{-1}\,d(h_\#\bfP_\theta)\\
			& =\int_X s_{\theta,V}s_{\theta,W}\,d\bfP_\theta\\
			& =\gfr(V,W)_\theta.
		\end{align*}
		Hence $\varphi_h^*(\gfr)=\gfr$. Since $\varphi_h$ is a $C^1$-diffeomorphism, it preserves the Riemannian distance. As $\gfr$ is smooth, the \textit{Myers-Steenrod theorem}, see, e.g., \cite[Chapter II]{Kobayashi1995}, implies that $\varphi_h$ is smooth. Therefore, $\varphi_h\in\Isom(M,\gfr)$.
	\end{proof}
	For the given $G$-invariant statistical model $(M,X,\bfP)$ with $G\leq\operatorname{Inv}(X)$, the analysis above implies that there exists a group homomorphism $\Phi:G\to\Isom(M,\gfr)$ where $\Isom(M,\gfr)$ is an isometry group on $(M,\gfr)$ equipped with the compact-open topology.
	The homomorphism is defined by $\Phi(h)=\varphi_h$ and 
	$$
	\Phi(h_1\circ h_2)=\Phi(h_1)\Phi(h_2)=\varphi_{h_1}\circ\varphi_{h_2}.
	$$
	In the remainder of the paper, we define $\calG:=\overline{\Phi(G)}\leq\Isom(M,\gfr)$ where $\overline{\Phi(G)}$ is the closure of $\Phi(G)$ in $\Isom(M,\gfr)$. We will see that the induced isometry subgroup $\calG$ can be preserved under equivariant and sufficient statistics.
	
	\begin{remark}
		Defining $\calG$ by $\overline{\Phi(G)}$ brings us many advantages, e.g., $\calG$ becomes a Lie group, which is a key ingredient in the argument for \cref{theorem_invariant_statistical_structures}. 
		However, a technical problem arises: for $\varphi\in\calG\setminus\Phi(G)$, if this $\varphi$ still preserves statistical structures ($\gfr$ and $\cac$)?
		We will see in \cref{subsection_quotient_Wasserstein_versus_Wasserstein_of_quotient} that, in general, the set of isometries preserving statistical structures is a proper subset of the whole isometry group.
		But for the induced group $\calG$, it always preserves the desired statistical structures.
	\end{remark}
	
	Now, we will show that for a $G$-invariant model $(M,X,\bfP)$, equivariant and sufficient statistics induce a new $G$-invariant model and the induced isometry subgroup is exactly the same as the original one. 
	To prove that, an important observation about a sufficient statistic $T: X\to X$ for the model $(M, X,\mathbf{P})$ is given below. 
	We recall that if $T$ is sufficient in the sense of \cref{definition_geometric_sufficiency}, then $\gfr=\gfr^T$ on $M$ where $\gfr^T$ is the induced Fisher-Rao metric by the model $(M,X,T_{\#}\bfP)$.
	
	\begin{lemma}\label{lemma_sufficiency_and_identity_map}
		Assume that $\mathbf P:M\to\mathcal M(X,\mathscr B)$ is $C^1$ and its scores are $L^2$-integrable. Let $\scrB_S\subseteq\scrB$ induced by $T$ be sufficient in the sense of \cref{definition_sufficient_sigma_algebra}. Then, for every $\theta\in M$ and $V\in T_\theta M$,
		\begin{equation}\label{equation_condition_expectation_of_sufficient_sub_algebra}
			s_{\theta,V}=\bbE_{\theta}(s_{\theta,V}\mid\scrB_S)\ \ \mathbf P_\theta\operatorname{-a.e.}
		\end{equation}
		Consequently, restriction to $\scrB_S$ preserves the Fisher-Rao metric. In particular, we have that $\gfr^{T}=\gfr$ where $\gfr^T$ is the Fisher-Rao metric on $M$ induced by the parameterized model $(M,X,T_{\#}\bfP)$.
	\end{lemma}
	
	\begin{proof}
		Fix $\theta\in M$ and $V\in T_\theta M$, for notational simplicity, we write
		$$
		s=s_{\theta,V},\ \ u=\bbE_\theta(s\mid\scrB_S).
		$$
		For every bounded $\scrB$-measurable function $f$, sufficiency provides a bounded $\scrB_S$-measurable function $\psi_f$ depending on $f$ but not on $\vartheta\in M$, such that
		$$
		\psi_f=\bbE_\vartheta(f\mid\scrB_S)\ \ \bfP_\vartheta\operatorname{-a.e.}
		$$
		In particular,
		$$
		\int_X f\,d\bfP_\vartheta=\int_X \psi_f\,d\bfP_\vartheta\ \ \operatorname{for\ ever\ }\vartheta\in M.
		$$
		
		Next, define the linear functional $F_f:\calM(X,\scrB)\to\bbR$ by $F_f(\mu):=\int_{X}fd\mu$ for the given Bounded $\scrB$-measurable $f$. Note that this functional is indeed a bounded linear functional with respect to the total variation norm in $\calM(X,\scrB)$, which implies it is $C^1$. We therefore define the function on $M$ by
		$$
		F_f\circ\bfP:M\to\bbR:\vartheta\mapsto\int_X fd\bfP_{\vartheta}.
		$$
		Since $\bfP$ is $C^1$, and $F_f$ is independent of $\vartheta$, we can calculate the differential of $F_f\circ\bfP$ at $\vartheta=\theta\in M$ by
		\begin{equation*}
			D_{\theta}(F_f\circ\bfP)(V)=F_f\circ D_{\theta}\bfP(V)=\int_{X}fd(D_{\theta}\bfP(V))=\int_{X}fsd\bfP_{\theta}.
		\end{equation*}
		for any $V\in T_{\theta}M$. Hence, differentiating $\int_X f\,d\mathbf P_\vartheta =\int_X \psi_f\,d\mathbf P_\vartheta$ at $\vartheta=\theta$ in direction $V$ gives the first equality below:
		$$
		\bbE_\theta(fs)=\bbE_\theta(\psi_fs)=\bbE_\theta(\psi_fu)=\bbE_\theta(fu).
		$$
		The second equality follows from the $\scrB_S$-measurability of $\psi_f$ and the definition $u=\bbE_\theta(s\mid\scrB_S)$. The third equality follows from the $\scrB_S$-measurability of $u$ and the identity $\psi_f=\bbE_\theta(f\mid\scrB_S)$. If $u$ is unbounded, the latter equality follows by truncation and approximation, since $u\in L^1(X,\bfP_\theta)$.
		
		Thus $\bbE_\theta(f(s-u))=0$ for every bounded measurable $f$. This yields $s=u$ $\bfP_\theta$-a.e., proving \eqref{equation_condition_expectation_of_sufficient_sub_algebra}. Finally, let $s_{\theta,V}^T$ be the score of $T_\#\bfP$. Since push-forward by the fixed statistic $T$ is a bounded linear operator on the space of finite signed measures $(\calM(X),\|\cdot\|_{\operatorname{TV}})$, the chain rule gives
		$$
		D_\theta(T_\#\bfP)(V)=T_\#(D_\theta\bfP(V)).
		$$
		Therefore, defining $\bfQ_{\theta}$ by $T_{\#}\bfP_{\theta}$, we have that $s^T:=s^T_{\theta,V}=\frac{d(D_{\theta}\bfQ(V))}{d(\bfQ(\theta))}$. By the identity above, we have $D_{\theta}\bfQ(V)=T_{\#}(D_{\theta}\bfP(V))=T_{\#}(s\bfP_{\theta})$. On the other hand, we have $D_{\theta}\bfQ(V)=s^T\bfQ_{\theta}$, this implies that for any $A\in\scrB$, we have that
		\begin{align*}
			D_{\theta}\bfQ(V)(A)= \int_{A}s^T_{\theta,V}d\bfQ_{\theta} = \int_{T^{-1}(A)} s_{\theta,V}^T\circ Td\bfP_{\theta} = \int_{T^{-1}(A)}s_{\theta,V}d\bfP_{\theta}
		\end{align*}
		where the last equality follows from the fact that $D_{\theta}\bfQ(V)(A)=T_{\#}(s\bfP_{\theta})(A)=s\bfP_{\theta}(T^{-1}(A))$.
		Since $s^T\circ T$ is $\scrB_S$-measurable and integrable, combining with the definitions of the score and conditional expectation and the uniqueness of the Radon-Nikodym derivative gives
		$$
		s_{\theta,V}^T\circ T=\bbE_\theta(s_{\theta,V}\mid\scrB_S)=s_{\theta,V}\ \ \bfP_\theta\operatorname{-a.e.}
		$$
		Therefore, for any $V,W\in T_\theta M$,
		$$
		\gfr^T(V,W)_\theta=\int_X (s_{\theta,V}^T\circ T)(s_{\theta,W}^T\circ T)\,d\bfP_\theta=\gfr(V,W)_\theta.
		$$
		This concludes the proof.
	\end{proof}
	
	\begin{remark}\label{remark_classic_induces_geometric_sufficiency}
		We give one observation on the preceding lemma.
		Under the regularity assumptions of \cref{lemma_sufficiency_and_identity_map}, classical sufficiency in the sense of \cref{definition_sufficient_sigma_algebra} implies preservation of Fisher information, i.e., sufficiency in the sense of \cref{definition_geometric_sufficiency} because
		$$
		\|s_{\theta,V}^T\|_{L^k(X,T_\#\bfP_\theta)}^k=\int_X |s_{\theta,V}^T\circ T|^k\,d\bfP_\theta=\|s_{\theta,V}\|_{L^k(X,\bfP_\theta)}^k.
		$$
	\end{remark}
	
	\begin{theorem}\label{theorem_geometric_analogy}
		Let $(M,X,\mathbf{P})$ and $G\leq\operatorname{Inv}(X)$ be a $G$-invariant statistical model. Suppose that $T: X\to X$ is an equivariant and sufficient statistic. Then the isometry $\varphi_h \in \operatorname{Isom}(M, \gfr)$ induced by the group action on the original model is an isometry for the model $(M, X, T_{\#}\mathbf{P})$. That is, for any $\theta\in M$, it holds $h_{\#}T_{\#}\bfP(\theta)=T_{\#}\bfP(\varphi_h(\theta))$.
	\end{theorem}
	\begin{proof}
		To prove the statement, we first claim that under our assumptions, for any fixed $h\in G\leq\operatorname{Inv}(X)$, there is $h_{\#}T_{\#}\mathbf{P}(M)=T_{\#}\mathbf{P}(M)$. This is straightforward by applying 
		$$
		h_{\#}T_{\#}=(h\circ T)_{\#}=(T\circ h)_{\#}=T_{\#}h_{\#}
		$$
		where the first and last equality follow from the \textit{fundamental identity} in optimal transport, see \cite[Lemma 1.2.7]{Figalli2023} for instance. The second equality follows from the $G$-equivariance of $T$. Thus combining with the assumption $h_{\#}\mathbf{P}(M)=\mathbf{P}(M)$, we have the desired equality $h_{\#}T_{\#}\mathbf{P}(M)=T_{\#}\mathbf{P}(M)$.
		
		Next, from \cref{proposition_isometry_group_for_FR_metric}, we know that the group action $h_{\#}$ on $\mathbf{P}(M)$ induces an isometric parameter transform $\varphi_h:M\to M$ such that for any $\theta\in M$, there is $h_{\#}\mathbf{P}(\theta)=\mathbf{P}(\varphi_h(\theta))$ and $\varphi_h\in\operatorname{Isom}(M,\gfr)$. Now we apply the push-forward via $T$ to both sides. Then we have
		$$
		h_{\#}T_{\#}\mathbf{P}(\theta)=T_{\#}h_{\#}\mathbf{P}(\theta)=T_{\#}\mathbf{P}(\varphi_h(\theta))
		$$
		where the first equality follows from $h_{\#}T_{\#}=T_{\#}h_{\#}$, the second equality follows form the fact that $h_{\#}\mathbf{P}(\theta)=\mathbf{P}(\varphi_h(\theta))$. This equality states that for any given $h\in G\leq\operatorname{Inv}(X)$, the group action $h_{\#}$ on the induced model $(M,X,T_{\#}\mathbf{P})$ yields the \textit{exact same} parameter transform $\varphi_h$ on $M$.
		
		Finally, according to \cref{lemma_sufficiency_and_identity_map}, since $\gfr^{T}=\gfr$, we conclude that $\operatorname{Isom}(M,\gfr^{T})=\operatorname{Isom}(M,\gfr)$. Hence, $\varphi_h\in\operatorname{Isom}(M,\gfr)$ implies that $\varphi_h\in\operatorname{Isom}(M,\gfr^T)$.
	\end{proof}
	
	Theorem \ref{theorem_geometric_analogy} provides a compatibility between the group invariance of a given model $(M,X,\bfP)$ and an equivariant and sufficient statistic $T:X\to X$.
	In our setting, an equivariant and sufficient statistic $T$ preserves both the Fisher-Rao metric and the induced parameter action. The compatibility is expressed by the following commutative diagram:
	\[
	\begin{tikzcd}[column sep=huge,row sep=large]
		\mathbf{P}(\theta)
		\arrow[r,"h_{\#}"]
		\arrow[d,"T_{\#}"']
		&
		\mathbf{P}(\varphi_h(\theta))
		\arrow[d,"T_{\#}"]
		\\
		T_{\#}\mathbf{P}(\theta)
		\arrow[r,"h_{\#}"']
		&
		T_{\#}\mathbf{P}(\varphi_h(\theta)).
	\end{tikzcd}
	\]
	The diagram follows from equivariance; sufficiency additionally ensures preservation of the Fisher-Rao metric. 
	This diagram concerns the commutativity between equivariant and sufficient statistics and the induced isometric group action, which is slightly different from the Hall-Wijsman-Ghosh's theorem.
	
	\subsection{Invariant statistical structures and homogeneity}\label{subsection_quotient_Wasserstein_versus_Wasserstein_of_quotient}
	
	In the previous subsection, we have seen that the $G$-invariance of a statistical model $(M,X,\bfP)$ induces an isometry subgroup on the associated statistical manifold.
	In this subsection, we show that the induced isometry subgroup preserves the statistical structures induced by the model.
	As a byproduct, we will also see that if the isometric action on the manifold induced by the group action on the sample space is transitive, i.e., $M\cong \calG/H$ where $H$ is the stabilizer of some $\theta_0\in M$ within $\calG$, then the statistical manifold is automatically homogeneous in the sense of \cite{Inoguchi2025}.
	
	Homogeneity is a key ingredient in studying symmetric properties of Riemannian manifolds. The homogeneous Riemannian manifolds have been well studied; see, e.g., \cite[Chapter 3]{Lee2018} or \cite[Chapter IV]{Helgason2001}.
	The homogeneous statistical manifold, which is a homogeneous Riemannian manifold and the statistical structures on it are invariant under the isometry subgroup action, has been explored by Inoguchi and Ohno \cite{Inoguchi2025} recently.
	We will show that, even though the action of the induced isometry subgroup $\calG:=\overline{\Phi(G)}$ does not need to be transitive on $M$, if the action is transitive, the statistical manifold $(M,\gfr,\cac)$ automatically becomes a homogeneous statistical manifold.
	
	Let us first briefly extend the results given in \cref{section_invariance_sufficiency_and_isometry} to other statistical structures.
	Recall that for a $k$-integrable model $(M,X,\bfP)$ with $k\geq3$, then Amari-Chentsov cubic tensor is given by
	$$
	\cac(U,V,W)_{\theta}:=\int_{X}s_{\theta,U}s_{\theta,V}s_{\theta,W}d\bfP_{\theta}(x)
	$$
	for any $\theta\in M$ and $U,V,W\in T_{\theta}M$. $s_{\theta,U}$ is the score function defined in \cref{section_preliminaries}.
	According to the argument in \cref{section_invariance_sufficiency_and_isometry}, it is not hard to observe that the induced isometry subgroup preserves $C$ at any $\theta\in M$ for any $U,V,W\in T_{\theta}M$, which can be formulated into the following theorem.
	\begin{theorem}\label{theorem_invariant_statistical_structures}
		Let $(M,X,\bfP)$ be a $G$-invariant model with $G\leq\operatorname{Inv}(X)$ having Borel measurable elements. Assume that $\calG:=\overline{\Phi(G)}\leq\Isom(M,\gfr)$. Then for every $\varphi\in\calG$, there is 
		\begin{equation}\label{equation_pullback_of_AC_tensor}
			\varphi^*\cac=\cac,
		\end{equation}
		i.e., the induced isometry subgroup preserves the Amari-Chentsov tensor.
	\end{theorem}
	\begin{proof}
		If $\varphi=\varphi_{h}=\Phi(h)$ for some $h\in G$, then repeating the argument in \cref{proposition_isometry_group_for_FR_metric} yields the desired statement.
		Therefore, it suffices to prove the statement for some $\varphi\in\calG\setminus\Phi(G)$.
		
		We recall that the manifolds that we use in the present paper are connected, complete and finite-dimensional. Hence, Myers-Steenrod's theory, see, e.g., \cite[Chapter 2.1]{Alexandrino2015}, says that $\calG\leq\Isom(M,\gfr)$ is a finite-dimensional Lie group whose manifold topology is the compact-open topology. Furthermore, the natural action 
		$$
		\zeta:\calG\times M\to M:(\varphi,\theta)\mapsto\zeta(\varphi,\theta):=\varphi(\theta)\in M
		$$
		is smooth. This yields that for an arbitrary sequence $\{\varphi_{n}\}\subset\Phi(G)$, which converges to the $\varphi\in\calG\setminus\Phi(G)$ with respect to the compact-open topology, the sequence of differentials $\{d\varphi_{n}\}$ converges to $d\varphi$ as well. 
		Therefore, for any $n\in\bbN_+$, $\varphi_n\in\Phi(G)$, we have
		$$
		(\varphi_{n})^*\cac=\cac\ \iff\ \cac\left( d\varphi_n(U),d\varphi_n(V),d\varphi_n(W) \right)_{\varphi_{n}(\theta)}=\cac(U,V,W)_{\theta}
		$$
		where $\theta\in M$ and $U,V,W\in T_{\theta}M$ are arbitrary.
		Hence, the fact that $\varphi_{n}\rightarrow\varphi$, $d\varphi_n\rightarrow d\varphi$ and $\cac$ is continuous yields
		$$
		\cac\left( d\varphi(U),d\varphi(V),d\varphi(W) \right)_{\varphi(\theta)}=\cac(U,V,W)_{\theta},
		$$
		which is equivalent to $\varphi^*\cac=\cac$ for $\varphi\in\calG$.
	\end{proof}
	This observation means that the isometry subgroup induced by the group action on the sample space preserves the statistical structures on the associated statistical manifold. Furthermore, \cref{theorem_geometric_analogy} implies that this property is stable under equivariant and sufficient statistics.
	The following simple example shows that the set of isometries that preserve $\gfr$ and $\cac$ in general is only a proper subset of the whole isometry group.
	
	\begin{example}\label{example_proper_subset_of_isometry_group}
		Let use use the \textit{exponential family} $p_\lambda(x)=\lambda e^{-\lambda x}$ on $(0,+\infty)$ as an example.
		By using the coordinate $t=\log\lambda$, we have that
		$$
		\gfr=dt^{\otimes2},\ \ \cac=-2dt^{\otimes3}
		$$
		Thus the reflection $r(t)=-t$ satisfies $r^*\gfr=\gfr$ but $r^*\cac=-\cac\neq\cac$. In contrast, the sample scalings $g_{a}(x)=ax$, $a>0$, induce translations $t\mapsto t-\log a$, which preserve both tensors.
	\end{example}
	
	As we can see, $\calG=\overline{\Phi(G)}\leq\Isom(M,\gfr)$ always preserve the statistical structures $(\gfr,\cac)$ induced by $(M,X,\bfP)$ on $M$.
	We can state this result in an alternative form by using $\alpha$-connections. For the notion of $\alpha$-connection $\nabla^{(\alpha)}$, see, e.g., \cite[Chapter 3.1]{Amari2000} or \cite[Chapter 2.5]{Ayetal2017}.
	\begin{proposition}\label{proposition_invariant_alpha_connection}
		Under the condition of \cref{theorem_invariant_statistical_structures}, the $\alpha$-connection $\AlphaC$ on the statistical manifold $(M,\gfr)$ associated to the model $(M,X,\bfP)$ is invariant under the action of $\calG$ for any $\alpha\in\bbR$. That is, $\varphi^*\AlphaC=\AlphaC$ for any $\alpha\in\bbR$, $\varphi\in\calG$. 
	\end{proposition}
	\begin{proof}
		We first note that since $\gfr$ is invariant under $\calG$, so the associated Levi-Civita connection $\LCcnt$ is also invariant under $\calG$, i.e., $\varphi^*\LCcnt=\LCcnt$ for any $\varphi\in\calG$.
		Then, combining with \cref{theorem_invariant_statistical_structures}, the relation between $\AlphaC$ and $\cac$ 
		\begin{equation*}
			\gfr(\AlphaC_UV,W)=\gfr(\LCcnt_UV,W)+\frac{\alpha}{2}\cac(U,V,W)
		\end{equation*}
		for any smooth vector field $U,V,W$ on $M$ implies that $\AlphaC$ is also invariant for any $\alpha\in\bbR$.
	\end{proof}
	
	These observation further implies that, if the induced isometry subgroup $\calG$ acts transitively on $M$, i.e., for any $\theta_1,\theta_2\in M$, there exists a $\varphi\in\calG$ such that $\theta_2=\varphi(\theta_1)$, then $(M,\gfr,\cac)$ is a homogeneous statistical manifold in the sense of \cite{Inoguchi2025}.

	\section{Reduction to the Group Invariant Statistical Model}\label{section_invariant_statistical_structures_and_homogeneity}
	
	In the previous section, we have seen that the $G$-invariance of a statistical model $(M,X,\bfP)$ induces an isometry subgroup on the associated statistical manifold. In this section, we will show that this $G$-invariance also induces a statistical reduction of the model $(M,X,\bfP)$.
	Then, we quantify the information loss caused by the reduction.
	
	To introduce the reduced model and the associated statistical model, let us first consider the restricted model.
	In this section, in addition to the smooth, isometric and proper action of $\calG$ on $M$, we further assume that $\calG$ \textit{freely} acts on $M$, i.e., $\Stab_{\calG}(\theta)=\{e\}$ for any $\theta\in M$. 
	Then, we can know that $N:=M/\calG$ is a smooth manifold, see, e.g., \cite[Chapter 21]{Lee2013} and \cite{Galaz2018}. Let $\gqout$ be the quotient metric on $N$ induced by $\gfr$ on $M$. Then $(N,\gqout)$ is a smooth Riemannian manifold of dimension $\operatorname{dim}M-\operatorname{dim}\calG$.
	
	To define the induced Fisher-Rao tensor on $N$, we need two new parameterized statistical models, the \textit{restricted model} $(M,X,\bfP^I)$ and the \textit{reduced model} $(N,X,\bfQ)$. 
	We also note that the restriction to $\scrB_I$ makes the Fisher-Rao tensor along $\calG$-orbit directions. The descended tensor on $N$ may still be degenerate, so its positive-definiteness will be assumed separately.
	We use $\phi:M\to N$ to represent the quotient map. Let $\scrB_I\subset\scrB$ be the invariant subalgebra of $\scrB$ with respect to $G$.
	We define the restricted map $\bfP^I:=\bfP|_{\scrB_I}$ by $\bfP^I(\theta)=\bfP(\theta)|_{\scrB_I}$. For any $A\in\scrB_I$, by the commutative property we saw in the previous subsection, we have
	\begin{align*}
		\bfP(\varphi_h(\theta))(A) = h_{\#}\bfP(\theta)(A) = \bfP(\theta)(h^{-1}(A)) = \bfP(\theta)(A) = \bfP(\theta)(A)|_{\scrB_I}
	\end{align*}
	for any $h\in G$. This implies that $\bfP^I:=\bfP|_{\scrB_I}$ is constant on each $\calG$-orbit in $M$. 
	This invariance extends from $\Phi(G)$ to its closure $\calG$. Indeed, if $\varphi_{h_j}\to f\in\mathcal G$, then total-variation continuity of $\bfP$ gives
	\begin{equation*}
		\bigl\|\bfP_{f(\theta)}|_{\scrB_I}-\bfP_\theta|_{\scrB_I}\bigr\|_{\mathrm{TV}}=\bigl\|\bfP_{f(\theta)}|_{\scrB_I}-\bfP_{\varphi_{h_j}(\theta)}|_{\scrB_I}\bigr\|_{\mathrm{TV}}\leq\bigl\|\bfP_{f(\theta)}-\bfP_{\varphi_{h_j}(\theta)}\bigr\|_{\mathrm{TV}}\longrightarrow 0.
	\end{equation*}
	Thus $\mathbf P^I$ is constant on $\mathcal G$-orbits and defines
	$$
	\bfQ:N\longrightarrow\calP(X,\scrB_I),\ \ \bfQ\circ\phi=\bfP^I.
	$$
	so that the parameterized statistical model $(N,X,\bfQ)$ is well-defined. 
	Hence $\mathbf Q$ is $C^1$ in the sense of \cite[Definition C.3]{Ayetal2017}.
	Therefore, we can define a new \textit{Fisher-Rao} tensor on $N$, which is denoted by $\gqfr$. Therefore, for any tangent vector $\upxi\in T_{\eta}N$, we define the \textit{score function} at $\eta\in N$ by 
	\begin{equation*}
		\bar{s}_{\eta,\upxi}=\frac{d(D_{\eta}\bfQ(\upxi))}{d\bfQ(\eta)}=:\partial_{\upxi}\log\bfQ(\eta).
	\end{equation*}
	Hence the induced Fisher-Rao tensor $\gqfr$ is defined by 
	\begin{equation*}
		\gqfr(\upxi_1,\upxi_2)_{\eta}:=\int_{X}\bar{s}_{\eta,\upxi_1}(x)\cdot \bar{s}_{\eta,\upxi_2}(x)d\bfQ_{\eta}(x)
	\end{equation*}
	for any $\upxi_1,\upxi_2\in T_{\eta}N$ and $Q_{\eta}:=Q(\eta)$.
	
	\begin{proposition}\label{proposition_conditional_expectation_of_quotient_score}
		For the parameterized statistical model $(M,X,\bfP)$ and $G\leq\operatorname{Inv}(X)$, we have
		\begin{equation}\label{equation_condition_expectation_of_quotient_score}
			s^I_{\theta,V}=\bbE_\theta(s_{\theta,V}\mid\scrB_I)\ \ \bfP_\theta^I\operatorname{-a.e.},
		\end{equation}
		where $s^I_{\theta,V}:=\frac{d(D_\theta\bfP^I(V))}{d\bfP^I(\theta)}$ is the score function of $\bfP^I$ with respect to $V\in T_\theta M$. Furthermore, $s^I_{\theta,V}\in L^2(X,\scrB_I,\bfP_\theta^I)$ for any $\theta\in M$ and $V\in T_\theta M$.
	\end{proposition}
	
	\begin{proof}
		We recall that $\bfP^I$ is a map from $M$ to $\calP(X,\scrB_I)$. Define the restriction operator
		$$
		\calR^I:\calM(X,\scrB)\to\calM(X,\scrB_I),\ \ \calR^I(\mu):=\mu|_{\scrB_I}.
		$$
		This operator is linear and bounded, since $\|\calR^I(\mu)\|_{\mathrm{TV}}\leq\|\mu\|_{\mathrm{TV}}$. Therefore, $\bfP^I=\calR^I\circ\bfP$ is still $C^1$ as a map into $\calM(X,\scrB_I)$, and
		$$
		D_\theta\bfP^I(V)=\calR^I\bigl(D_\theta\bfP(V)\bigr)=\bigl(D_\theta\bfP(V)\bigr)|_{\scrB_I}.
		$$
		
		For any $A\in\scrB_I$, by the definition of the score function and conditional expectation, we have
		\begin{align*}
			\bigl(D_\theta\bfP^I(V)\bigr)(A)
			=\int_A s_{\theta,V}\,d\bfP_\theta=\int_A\bbE_\theta(s_{\theta,V}\mid\scrB_I)\,d\bfP_\theta=\int_A\bbE_\theta(s_{\theta,V}\mid\scrB_I)\,d\bfP_\theta^I.
		\end{align*}
		The last equality follows because the integrand is $\scrB_I$-measurable and $\bfP_\theta^I=\bfP_\theta|_{\scrB_I}$.
		Since $A\in\scrB_I$ is arbitrary, the uniqueness of the Radon-Nikodym derivative gives
		$$
		s^I_{\theta,V}=\frac{d(D_\theta\bfP^I(V))}{d\bfP^I(\theta)}=\bbE_\theta(s_{\theta,V}\mid\scrB_I)\ \ \bfP_\theta^I\operatorname{-a.e.}
		$$
		This proves \eqref{equation_condition_expectation_of_quotient_score}.
		
		Finally, since conditional expectation is an orthogonal projection in $L^2(X,\scrB,\bfP_\theta)$,
		$$
		\|s^I_{\theta,V}\|_{L^2(X,\scrB_I,\bfP_\theta^I)}=\|\bbE_\theta(s_{\theta,V}\mid\scrB_I)\|_{L^2(X,\scrB,\bfP_\theta)}\leq\|s_{\theta,V}\|_{L^2(X,\scrB,\bfP_\theta)}.
		$$
		Hence $s^I_{\theta,V}$ is in $L^2(X,,\scrB_I,\bfP^I_{\theta})$ as desired.
	\end{proof}
	
	The next result reveals the relation between $s^I_{\theta,V}$ and $\bar{s}_{\eta,\upxi}$.
	
	\begin{proposition}\label{proposition_restricted_score_and_induced_score}
		Let $(M,X,\bfP)$ be the parameterized statistical model, and $(N,X,\bfQ)$ the reduced model on $(X,\scrB_I)$ as above, with $\bfP^I=\bfQ\circ\phi$ and $\phi:M\to N$ a smooth submersion. For any $\eta\in N$, $\upxi\in T_\eta N$, and $\theta\in\phi^{-1}(\eta)$, let $V\in T_\theta M$ be the horizontal lift of $\upxi$. Then
		\begin{equation}\label{equation_restricted_score_and_induced_score}
			\bar{s}_{\eta,\upxi}=s^I_{\theta,V}=\bbE_\theta(s_{\theta,V}\mid\scrB_I)\ \ \bfQ(\eta)\operatorname{-a.e.}
		\end{equation}
	\end{proposition}
	
	\begin{proof}
		We first recall that $\bfP^I=\bfQ\circ\phi$. Since $\phi$ is a smooth submersion, according to \cite[Lemma 1.7]{Glockner2015} or \cite[Exercise 1.7.5]{Schmeding2023} which is indeed \cite[Proposition 4.1.13]{Margalef1992}, for any $\eta\in N$, there is an open neighborhood $\eta\in U\subset N$ and an associated smooth map $j:U\to M$ such that $\phi\circ j=\Id_N|_U$. This identity implies that for any $\eta\in N$, there is
		$$
		\bfP^I(j(\eta))=(\bfQ\circ\phi\circ j)(\eta)=\bfQ(\eta)\ \ \implies\ \ \bfQ|_U=\bfP^I\circ j.
		$$
		Thus, according to \cite[Lemma 1.8]{Glockner2015} or \cite[Exercise 1.7.6]{Schmeding2023}, $\bfQ$ is $C^1$ as a map into $\calM(X,\scrB_I)$. 
		
		Fix $\eta\in N$, $\upxi\in T_\eta N$, and $\theta\in\phi^{-1}(\eta)$. Let $V\in T_\theta M$ be the horizontal lift of $\upxi$, so that $D_\theta\phi(V)=\upxi$. By the chain rule, we have
		$$
		D_\theta\bfP^I(V)=D_{\phi(\theta)}\bfQ\bigl(D_\theta\phi(V)\bigr)=D_\eta\bfQ(\upxi).
		$$
		Moreover, $\bfP^I(\theta)=\bfQ(\eta)$. Therefore, by the definition of the score function and \cref{proposition_conditional_expectation_of_quotient_score},
		$$
		\bar{s}_{\eta,\upxi}=\frac{d(D_\eta\bfQ(\upxi))}{d\bfQ(\eta)}=\frac{d(D_\theta\bfP^I(V))}{d\bfP^I(\theta)}=s^I_{\theta,V}=\bbE_\theta(s_{\theta,V}\mid\scrB_I)\qquad\bfQ(\eta)\operatorname{-a.e.}
		$$
		This proves \eqref{equation_restricted_score_and_induced_score}.
	\end{proof}
	
	A useful observation of the analysis above can be formulated as follows.
	\begin{corollary}\label{corollary_relation_between_quotient_and_restricted_FR}
		For the parameterized statistical model $(M,X,\bfP)$, denote the restricted model by $(M,X,\bfP^I)$ and the associated Fisher-Rao tensor by $\gfr^I(V_1,V_2):=\bbE_{\theta}(s^I_{\theta,V_1}s^I_{\theta,V_2})$. Then we have
		\begin{equation*}
			(\phi^*\gqfr)(V_1,V_2)_{\theta}=\gqfr(D_{\theta}\phi(V_1),D_{\theta}\phi(V_2))_{\phi(\theta)}=\gfr^I(V_1,V_2)_{\theta}
		\end{equation*}
	\end{corollary}

	\begin{proposition}\label{proposition_Hall_invariant_Fisher_tensor}
		Under the conditions of \cref{ref.Hall_theorem_3.1}, define $\bfP^{SI}_\theta:=\bfP_\theta|_{\scrB_{SI}}$, and let $\gfr^{SI}$ denote its Fisher-Rao tensor.
		Then
		$$
		\gfr^{SI}=\gfr^I=\phi^*\gqfr.
		$$
		More precisely, for every $\theta\in M$ and $V\in T_\theta M$,
		$$
		\bbE_\theta(s_{\theta,V}\mid\mathscr B_{SI})=\bbE_\theta(s_{\theta,V}\mid\mathscr B_I)\ \ \mathbf P_\theta\operatorname{-a.e.}
		$$
	\end{proposition}
	\begin{proof}
		By Hall-Wijsman-Ghosh's theorem, cf. \cref{ref.Hall_theorem_3.1}, $\mathscr B_{SI}$ is sufficient for $\mathscr B_I$.
		Applying the arguments in \cref{lemma_sufficiency_and_identity_map} to this restricted experiment gives
		$$
		s^I_{\theta,V}=\bbE_\theta(s^I_{\theta,V}\mid\scrB_{SI}).
		$$
		Since $\scrB_{SI}\subseteq\scrB_I$ is sufficient, by \eqref{equation_condition_expectation_of_quotient_score} and \eqref{equation_condition_expectation_of_sufficient_sub_algebra}, we have the restricted-score identities $s_{\theta,V}^I=\bbE_{\theta}(s_{\theta,V}|\scrB_I)$, and $s_{\theta,V}^{SI}=\bbE_{\theta}(s_{\theta,V}^I|\scrB_{SI})$, respectively. 
		Then the basic property of conditional expectation $\bbE_{\theta}(\bbE_{\theta}(\cdot|\scrB)|\scrF)=\bbE_{\theta}(\cdot|\scrF)$ if $\scrF\subset\scrB$ shows that
		$$
		\bbE_\theta(s^I_{\theta,V}\mid\scrB_{SI})=\bbE_{\theta}(\bbE_{\theta}(s_{\theta,V}|\scrB_I)|\scrB_{SI})=\bbE_\theta(s_{\theta,V}\mid\scrB_{SI}).
		$$
		Thus the two restricted scores agree, and consequently $\gfr^{SI}=\gfr^I$. The remaining equality is the pullback identity, which follows from \cref{corollary_relation_between_quotient_and_restricted_FR}.
	\end{proof}
	\begin{remark}
		Under the conditions of the Hall-Wijsman-Ghosh's theorem, $\scrB_S\cap\scrB_I$ is sufficient for the invariantly restricted model $\bfP^I$.
		Consequently, the preceding result yields $g_{SI}=g_I$: the further sufficient restriction introduces no additional loss of Fisher information.
		Thus, for the Fisher-geometric analysis below, we continue to work with $(M,X,\bfP^I)$ rather than $(M,X,\bfP^{SI})$.
	\end{remark}
	
	In the remainder of the section, we assume that $\gqfr$ is a smooth metric on $N$.
	The following result measures the loss of Fisher information under restriction to $\scrB_I$, evaluated along the $\gfr$-horizontal lift of a quotient tangent vector.
	
	\begin{theorem}\label{theorem_comparison_of_quotient_and_induced_FR}
		Let $N:=M/\calG$ be the smooth quotient manifold as above, and let $\gqout$ and $\gqfr$ be the quotient Riemannian metric and induced Fisher-Rao metric as above, respectively. Then we have
		\begin{equation}\label{equation_comparison_of_quotient_and_induced_FR}
			\gqout(\upxi,\upxi)_{\eta}-\gqfr(\upxi,\upxi)_{\eta}=\bbE_{\theta}(\Var(s_{\theta,V}|\scrB_I))
		\end{equation}
		where $\bbE_{\theta}(\Var(s_{\theta,V}|\scrB_I)):=\left\| s_{\theta,V}-\bbE_{\theta}(s_{\theta,V}|\scrB_I) \right\|_2^2$, $\theta\in\phi^{-1}(\eta)$ and $V\in T_{\theta}M$ is the horizontal lift of $\upxi$. 
	\end{theorem}
	\begin{proof}
		The fact that $\bbE(\cdot|\scrB_I)$ is an orthogonal projection in $L^2(X,\scrB,\bfP_{\theta})$ implies that for any $\upxi\in T_{\eta}N$, the following inequality naturally holds:
		\begin{equation*}
			\gqfr(\upxi,\upxi)_{\eta}=\left\| \bbE_{\theta}(s_{\theta,V}|\scrB_I) \right\|_{L^2(X,\scrB,\bfP_{\theta})}^2\leq\left\| s_{\theta,V} \right\|_{L^2(X,\scrB,\bfP_{\theta})}^{2}=\gqout(\upxi,\upxi)_{\eta}.
		\end{equation*}
		Where $\theta\in\phi^{-1}(\eta)$ and $V\in T_{\theta}M$ is the horizontal lift of $\upxi\in T_{\eta}N$. The first equality follows from \cref{proposition_conditional_expectation_of_quotient_score,proposition_restricted_score_and_induced_score}; the second equality follows from the definition of the quotient Riemannian metric.
		
		Therefore, we note that by the orthogonality of $\bbE_{\theta}(\cdot|\scrB_I)$, we have that
		\begin{align*}
			\left\| s_{\theta,V}-\bbE_{\theta}(s_{\theta,V}|\scrB_I) \right\|_2^2 & = \langle s_{\theta,V}-\bbE_{\theta}(s_{\theta,V}|\scrB_I)\,,\,s_{\theta,V}-\bbE_{\theta}(s_{\theta,V}|\scrB_I)\rangle\\
			& = \langle s_{\theta,V}-\bbE_{\theta}(s_{\theta,V}|\scrB_I)\,,\,s_{\theta,V}\rangle\\
			& = \langle s_{\theta,V}\,,\,s_{\theta,V}\rangle-\langle \bbE_{\theta}(s_{\theta,V}|\scrB_I)\,,\,\bbE_{\theta}(s_{\theta,V}|\scrB_I)\rangle\\
			& = \gqout(\upxi,\upxi)_{\eta}-\gqfr(\upxi,\upxi)_{\eta},
		\end{align*}
		which implies the desired identity \eqref{equation_comparison_of_quotient_and_induced_FR}.
	\end{proof}
	
	We conclude the section by estimating the total information loss caused by the reduction, i.e., the difference between $\gfr$ and $\phi^{\ast}\gqfr$ on $M$.
	
	\begin{theorem}\label{theorem_estimation_of_information_loss_caused_by_restriction}
		Let $\gqfr$ be the induced Fisher-Rao metric on $N$ induced by the model $(N,X,\bfQ)$. Then we have
		\begin{equation}\label{equation_estimation_of_information_loss_caused_by_restriction}
			\gfr(V,V)_{\theta}-(\phi^{\ast}\gqfr)(V,V)_{\theta}=\bbE_{\theta}(\Var(s_{\theta,V}|\scrB_{I}))
		\end{equation}
		where $\theta\in M$, $V\in T_{\theta}M$ are arbitrary, $\phi:M\to N$ is the quotient map.
	\end{theorem}
	\begin{proof}
		According to \cref{proposition_conditional_expectation_of_quotient_score,corollary_relation_between_quotient_and_restricted_FR}, we can conclude that
		$$
		\gfr(V,V)_{\theta}=\bbE_{\theta}(s_{\theta,V}^2)=||s_{\theta,V}||^2_{L^2(X,\bfP_{\theta})}
		$$
		and
		$$
		(\phi^{\ast}\gqfr)(V,V)_{\theta}=\gfr^I(V,V)_{\theta}=\bbE_{\theta}\left(\left( \bbE_{\theta}(s_{\theta,V}|\scrB_I) \right)^2\right)=\| \bbE_{\theta}(s_{\theta,V}|\scrB_I) \|^2_{L^2(X,\bfP_{\theta})}.
		$$
		Hence, the computation in \cref{proposition_conditional_expectation_of_quotient_score,theorem_comparison_of_quotient_and_induced_FR} shows that
		$$
		\gfr(V,V)_{\theta}-(\phi^{\ast}\gqfr)(V,V)_{\theta}=\left\| s_{\theta,V}-\bbE_{\theta}(s_{\theta,V}|\scrB_I) \right\|_2^2=\bbE_{\theta}(\Var(s_{\theta,V}|\scrB_{I})).
		$$
		Thus we have the desired identity \eqref{equation_estimation_of_information_loss_caused_by_restriction}.
	\end{proof}
	
	One can further observe that for the variance of the invariant observation appeared above, by some straightforward computation, for any $\eta\in N$, $\xi\in T_{\eta}N$ with $\theta\in\phi^{-1}(\eta)$ and the tangent vector $V\in T_{\theta}M$ such that $D_{\theta}\phi(V)=\xi$, it can be decomposed as follows
	\begin{equation*}\label{equation_decomposition_of_information_loss}
		\bbE_{\theta}(\Var(s_{\theta,V}|\scrB_{I}))=\gfr(V^{\operatorname{ver}},V^{\operatorname{ver}})_{\theta}+(\gqout-\gqfr)(\xi,\xi)_{\eta},
	\end{equation*}
	where $V=V^{\operatorname{hor}}+V^{\operatorname{ver}}$ and $V^{\operatorname{hor}}$, $V^{\operatorname{ver}}$ mean the horizontal and vertical parts of $V$ with respect to $\gfr$ in $T_{\theta}M$, respectively.

	\section{Application in circular data analysis}\label{section_application_in_circular_data_analysis}
	
	In this section, let us see an important example from \textit{circular data analysis} in statistics; see, e.g., \cite[Chapter 3]{Kanti2000} for details. We will see that for a parameterized statistical model $(M,X,\bfP)$ satisfying certain regularity conditions, our analysis given in \cref{section_invariance_sufficiency_and_isometry} can provide a quantitative description of information loss caused by the reduction.
	
	We first construct a concrete example based on the \textit{generalized cardioid distribution} studied in \cite{Paual2021}.
	Then we will give the information loss caused by the reduction induced by the model's $G$-invariance.
	
	In this section, we set $M:=\bbT^2$ where $\bbT:=\bbR/2\pi\bbZ$ and $X=\bbT^2$. Let $\bfX:=(\bfX_1,\bfX_2)$ be a random variable in $X$ such that $\bfX_1$ and $\bfX_2$ are independent. For any given $\theta=(\theta_1,\theta_2)$, we denote the independent random variable by $\bfX_1$ and $\bfX_2$ with density 
	$$
	p_{\theta_i}(x_i):=\frac{1}{2\pi}\left( 1+a\cos(x_i-\theta_i) \right),\ 0<a<1,\ i=1,2.
	$$
	We note that, for simplicity, we consider this as the density function of $\bfX_i$ with respect to the standard Lebesgue measure on $[0,2\pi)$, which can be viewed as the push-forward of the Riemannian volume measure on $\bbT$ via the angular coordinate $x\in[0,2\pi)$. Therefore, by the independence of $\bfX_1$ and $\bfX_2$ implies that
	$$
	\bfP(\theta)=p_{\theta}(x)dx\in\calP(X),\ \operatorname{and}\ p_{\theta}(x)=p_{\theta_1}(x_1)p_{\theta_2}(x_2).
	$$
	In this setting, although $M=X=\bbT^2$, they have different meanings: $\theta\in \bbT^2$ means the parameter of a certain cardioid distribution, while $\bfX\in \bbT^2$ means that the random variable under consideration represents a random angle.
	
	We select the group action $G\leq\operatorname{Inv}(X)$ such that $G=\bbT$, i.e., the \textit{rotation action} on $\bbT^2$. That is, for any $\alpha\in\bbT$, define $h_{\alpha}\in G$ and acts on $X$ by $h_{\alpha}(x_1,x_2)=(x_1+\alpha,x_2+\alpha)$. It is straightforward to verify that this $(M,X,\bfP)$ and $G\leq\operatorname{Inv}(X)$ in this setting satisfy the assumptions in \cref{section_invariance_sufficiency_and_isometry}, so we omit the verification here.
	
	By the definition of $(M,X,\bfP)$ and some elementary computation, one can observe that the induced isometry subgroup $\calG$ of $\Isom(M)$ is itself, as we have seen in \cref{section_invariance_sufficiency_and_isometry}, there is
	$$
	(h_{\alpha})_{\#}\bfP(\theta_1,\theta_2)=\bfP(\varphi_{\alpha}(\theta_1),\varphi_{\alpha}(\theta_2))=\bfP(\theta_1+\alpha,\theta_2+\alpha).
	$$
	Moreover, this action $\calG$ on $M$ is \textit{free, isometric and smooth}, which implies $N:=M/\calG$ is also a smooth Riemannian manifold.
	
	We use $\theta$ and $\bfX$ to represent $(\theta_1,\theta_2)$ and $(\bfX_1,\bfX_2)$, respectively. Since we have assumed $\bfX_1$ and $\bfX_2$ are independent, we can denote the joint density function of $\bfX$ by $p_{\theta}(x)=p_{\theta_1}(x_1)p_{\theta_2}(x_2)$. Let us compute the Fisher-Rao metric on $M$. Since $\bfX_1$ and $\bfX_2$ are one-dimensional random variable determined by one-dimensional parameter $\theta_1$ and $\theta_2$, we have that for $V=\partial_{\vartheta}$, there is
	$$
	s_{\vartheta,V}(x)=\partial_V\log(\bfP(\vartheta))=\frac{a\sin(x-\vartheta)}{1+a\cos(x-\vartheta)}.
	$$
	Therefore, the Fisher information, i.e., the entries of the Fisher-Rao metric, is given by
	$$
	I(a)=\frac{1}{2\pi}\int_{0}^{2\pi}\frac{a^2\sin[2](u)}{1+a\cos(u)}du=1-\sqrt{1-a^2}.
	$$
	The notation $du$ here can be understood as the standard Lebesgue measure on $[0,2\pi)$, but it can be viewed as the Riemannian volume measure $\Vol$ on $\bbT$ via the angular coordinate as well. 
	To show the application of \cref{theorem_comparison_of_quotient_and_induced_FR,theorem_estimation_of_information_loss_caused_by_restriction}, we select a certain parameter $a:=\frac{\sqrt{3}}{2}$, which implies $I(a)=\frac{1}{2}$.
	
	Since $\bfX_1$ and $\bfX_2$ are independent and the integral of the score function is zero, we have
	$$
	\bbE_{\theta}\left( s_{\theta,\partial_{\theta_1}}(\bfX_1)s_{\theta,\partial_{\theta_2}}(\bfX_2) \right)=0\ \implies\ \gfr=\frac{1}{2}(d\theta_1^2+d\theta_2^2).
	$$
	This implies that $\varphi_{\alpha}$ given above is indeed an isometry for any $\alpha\in\bbT$. This is also a flat metric, which is essentially the same as the usual Riemannian flat metric on the torus, and we have
	$$
	\sfd_{\gfr}(\theta,\vartheta)=\frac{\sqrt{2}}{2}\sqrt{\sfd_{\bbT}(\theta_1,\vartheta_1)^2+\sfd_{\bbT}(\theta_2,\vartheta_2)^2}
	$$
	where $\sfd_{\bbT}(x,y)=\min_{k\in\bbZ}|x-y+2k\pi|$ is the Riemannian distance on $\bbT$ induced by the restricted Euclidean metric.
	
	Next, we define a new random variable $\bfD:=\bfX_2-\bfX_1\ \operatorname{mod}(2\pi)$, as a sample, it means $\bfD(x_1,x_2)=x_2-x_1\ \operatorname{mod}(2\pi)$. The statistic $\bfD(x_1,x_2)$ is a maximum invariant under simultaneous rotations.
	Therefore, we can define the \textit{invariant subalgebra} of $X$ with respect to the group action $G$ by $\scrB_I=\sigma(D)$. So, we observe that $\bfD$ realizes the invariant experiment in this example.
	Hence, we can observe that $N:=M/\calG\cong\bbT$ and the quotient map $\phi$ is given by $\phi(\theta)=\eta:=\theta_2-\theta_1$.
	
	According to some elementary calculus in differential geometry, for a tangent vector $\xi=v\partial_{\eta}$ on $N$, its horizontal lift on $M$ is given by $V:=\frac{v}{2}(-\partial_{\theta_1}+\partial_{\theta_2})$, which implies the quotient Riemannian metric on $N$ induced by $\gfr$ on $M$ is given by  $\gqout=\frac{1}{4}d\eta^2$.
	On the other hand, the density function of $\bfD$ can be calculated by convolution as below
	$$
	p_{\eta}(d)=\int_{0}^{2\pi}p_{\theta_1}(x)p_{\theta_2}(x+d)dx=\frac{1}{2\pi}\left( 1+\frac{3}{8}\cos(d-\eta) \right).
	$$
	Then we compute the Fisher information, which is the Fisher-Rao metric in the one-dimensional case, to obtain $\gqfr=\left( 1-\frac{\sqrt{55}}{8} \right)d\eta^2$.
	This shows a strict inequality $0<\gqfr<\gqout$. Therefore, by applying \cref{theorem_comparison_of_quotient_and_induced_FR}, for any given $\xi=v\partial_{\eta}\in T_{\eta}N$ with $v\neq0$, there is
	$$
	\bbE_{\theta}(\Var(s_{\theta,V}|\scrB_I))=\gqout(\xi,\xi)-\gqfr(\xi,\xi)=\left( \frac{\sqrt{55}-6}{8} \right)v^2>0.
	$$
	For a general tangent vector $V\in T_{\theta}M$ of the form $V=a\partial_{\theta_1}+b\partial_{\theta_2}$, one can observe that
	$$
	V^{\operatorname{ver}}=\frac{a+b}{2}(\partial_{\theta_1}+\partial_{\theta_2}),\ V^{\operatorname{hor}}=\frac{b-a}{2}(-\partial_{\theta_1}+\partial_{\theta_2}),
	$$
	which implies that the total information loss is 
	$$
	\calL_{\theta}(V)=\frac{(a+b)^2}{4}+\frac{\sqrt{55}-6}{8}(b-a)^2.
	$$
	Therefore, in our example, $a=-v/2$, $b=v/2$, the total information loss in \cref{theorem_estimation_of_information_loss_caused_by_restriction} is the same as the information loss in the horizontal direction.
	
	To understand the meaning of the information loss, we can consider the statistical meaning of the random variable $\bfD$ as follows.
	The statistic \(\bfD=\bfX_2-\bfX_1\) is a maximal invariant under simultaneous rotations, and its distribution remains identifiable in the relative parameter \(\eta=\theta_2-\theta_1\). 
	Nevertheless, \(g_{\mathrm{inv}}<g_{\mathrm{quot}}\) shows that retaining all invariant observations does not preserve the Fisher information carried by horizontal parameter variations. 
	The conditional variance of the score quantifies this loss. In this example, the invariant experiment retains the constant fraction \(4-\sqrt{55}/2\) of the horizontal Fisher information. 
	
	From an inferential perspective, the reduced model provides a likelihood for the relative parameter $\eta$ that does not depend on the common angular origin.
	It therefore supports estimation and hypothesis testing for relative orientation when absolute orientation is unavailable or intentionally disregarded. 
	The quantity $\mathbb{E}_{\theta}(\operatorname{Var}(s_{\theta,V}\mid\mathscr{B}_{I}))$, evaluated along horizontal directions, measures the local information cost of this restriction: it records the component of the score that cannot be recovered from the invariant observation. Nevertheless, the positive Fisher information $I_{\mathrm{inv}}=1-\sqrt{55}/8$ ensures that the reduced model retains sensitivity to changes in $\eta$. 
	For $n$ independent repetitions in which only the differences $(\bfD_1,\ldots,\bfD_n)$ are retained, standard likelihood theory yields, in a local angular coordinate, asymptotically efficient estimation with variance scale $(nI_{\mathrm{inv}})^{-1}$. 
	Thus, the reduction preserves the possibility of consistent inference while reducing the precision available from each observation. This example consequently provides a tractable benchmark for evaluating the statistical cost of using relative measurements, and for studying estimation and testing procedures that are unaffected by a common rotation within each observed pair.

	
	\section*{AI-declaration}
	The AI tool (GPT-6 Astra in the ChatGPT Plus plan) is used to check the English text, the consistency and completeness of the paper, and to search for related references.
	All AI-generated suggestions were independently verified by the author, who takes full responsibility for the content and accuracy of the final manuscript.
	
	\section*{Acknowledgments}\label{section_acknowledgements}
	The author thanks his PhD supervisor, Prof. Hiroshi Hirai, for helpful comments and suggestions.
	This work is supported by JST SPRING, Grant Number JPMJSP2125. 
	
	\printbibliography 
	
\end{document}